\input amstex
\documentstyle {amsppt}
\pageheight{50.5pc} \pagewidth{32pc} %
\topmatter
\title{On Markov property of strong solutions to SDE with generalized coefficients}
\endtitle
\author
 Ludmila L. Zaitseva
\endauthor
\address
Kyiv Taras Shevchenko university, Volodymyrska 64, Kyiv, Ukraine,
01033
\endaddress
\abstract We show the complete proof of the Markov property of the
strong solution to a multidimensional SDE whose coefficients
involve local time on a hyperplane of the unknown process.
\endabstract
\email zaitseva\@univ.kiev.ua
\endemail
\endtopmatter
\document
\rightheadtext{Markov property of solutions to SDE with
generalized coefficients} \leftheadtext{Zaitseva L.L.}

\head
{1. Introduction.}
\endhead

We consider one class of stochastic differential equations with
the  local time on a hyperplane of the unknown process included
into the drift and martingale part. It is shown that this equation
has the strong unique solution and this solution is a Markov
process. The processes constructed may be described as a
generalized diffusion process in the sense of Portenko (see, for
example, \cite{1} ) with generalized both the drift vector and
diffusion matrix.

One of the most known example of a process described by an SDE
with its coefficients containing the local time of the unknown
process is the skew Brownian motion, introduced in \cite{2}
(Section 4.2, Problem 1). In 1981 J.M.~Harrisson, L.A.~Shepp (see
\cite{3}) proved that skew Brownian motion may be constructed as
the strong solution to the SDE of the form $dx(t)=qd\eta_t+dw(t),$
where $q\in[-1,1]$ is a given parameter and $\{\eta_t\}$ is the
local time at $0$ of the process $\{x(t)\}$. One-dimensional
equations with the drift of more general form were considered by
J.F.~Le~Gall, M.~Barlow, K.~Burdzy, H.~Kaspi, A.~Mandelbaum (see
\cite{4,5}). The main result of those papers is the theorem of the
existence and uniqueness of the strong solution to the
corresponding SDE.

The characterization of the obtained strong solution (say, as a
Markov process or generalized diffusion process) is a separate
complicated problem. J.M.~Harrisson and L.A.~Shepp in their paper
while stating that the skew Brownian motion, obtained as the
strong solution to the SDE, is a Markov process, referred to the
corresponding results for the standard SDE in \cite{6}. However,
this question is rather delicate, let us discuss it in a more
details. The scheme of proof in \cite{6} (or, in more generality,
in \cite{7}, Chapter 6) consists of two steps. The first step
contains the construction of a modification of the strong
solution, jointly measurable w.r.t. starting point and the Wiener
noise. The second one, namely the proof of the Markov property,
essentially use the fact that the increments of the Wiener process
are independent. The first step requires the proof of the
additional fact that the solution is continuous in probability as
a function of starting point. This proof is non-trivial  even in
the simplest case of a skew Brownian motion (see \cite{8} for a
general result) and at this point the arguments of J.M.~Harrisson,
L.A.~Shepp are incomplete. Notice that A.M.~Kulik (see \cite{9})
gives sketch of the proof of the Markov property in the case of
one-dimensional equation same to considered by J.F.~Le~Gall, which
takes this difficulty into account. In the multidimensional case
another difficulty appears, caused by the fact that the martingale
noise does  not necessarily have independent increments. The main
purpose of this paper is to give complete proof of the Markov
property for the strong solution to an SDE involving local time of
the unknown process in the most general  multidimensional case.

\head {2. The construction of the process.}
\endhead

Let $S$ be a hyperplane in $\Re^d$ orthogonal to the fixed unit
vector $\nu\in\Re^d:$ $S=\left\{x\in\Re^d|(x,\nu)=0\right\}.$ We
denote by $\pi_S$ the operator of orthogonal projection on $S$ and
by $L$ the one-dimensional subspace of $\Re^d,$ generated by
$\nu.$

We consider a Wiener process $\{w(t)\}$ in $\Re^d$ and filtration
${\Cal F}_t^w=\sigma \left\{w(u),0\leq u\leq t\right\}, t\geq 0.$
We denote $w(t)=(w^1(t),w^S(t)),$ where $w^1(t)=(w(t),\nu),
w^S(t)=\pi_S w(t).$ For a given parameter $q\in [-1,1]$ and
initial point $x^1_0\in L$ we construct skew Brownian motion (see
\cite{3}), i.e. a pair of $\{{\Cal F}_t^w\}$-adapted processes
$\{(x^1(t),\eta_t)\},$ where $\{\eta_t\}$ is the local time at $0$
of $\{x_1(t)\}:$
$$\eta_t=\lim_{\varepsilon\downarrow
0}{1\over2\varepsilon}\int_0^t{\hbox{\rm
1\!I}}_{\{|x^1(\tau)|\leq\varepsilon\}}d\tau,
$$
such that $\{x^1(t)\}$ and $\{\eta_t\}$ satisfy the equality
$$
x^1(t)=x^1_0+q\eta_t+w^1(t)
$$
for all $t\geq 0.$

Let us assume that there exists one more Wiener process
$\{\widetilde{w}(t)\}$ in $S$ that does not depend on $\{w(t)\}.$
Consider new process $\zeta(t)=\widetilde{w}(\eta_t), t\geq 0.$ We
shall deal with the SDE involving processes $\{w(t)\}$ and
$\{\zeta(t)\}$ as a martingale part of the equation. Therefore we
have to construct a filtration such that $\{\zeta(t)\}$ is the
square integrable martingale and $\{w(t)\}$ is the Wiener process
with respect to this filtration. Remark that the process
$\{\eta_t\}$ is not a stopping time w.r.t. $\{{\Cal F}_t^w\}$ and
therefore the classic results cannot be applied here

For $t\geq 0$ we consider
$$
\widetilde{{\Cal
F}}_t=\sigma\left\{\left\{\widetilde{w}(s)\in\Gamma\right\}\bigcap\left\{\eta_t\geq
s\right\},s\geq 0,\Gamma\in{\Cal B}_S\right\},
$$
where ${\Cal B}_S$ is the Borel $\sigma$-algebra on $S.$ We put
$${\Cal M}_t={\Cal F}^w_t\bigvee\widetilde{{\Cal F}}_t.$$

The following three lemmas shows that $\{{\Cal M}_t\}$ is the
required filtration.

\proclaim{Lemma 1}$\{{\Cal M}_t\}$ is a filtration.\endproclaim
\demo{Proof} We have to prove that ${\Cal M}_{t_1}\subseteq {\Cal
M}_{t_2},$ when $t_1\leq t_2.$ For fixed $s\geq 0$ and
$\Gamma\in{\Cal B}_S$ let us consider the set
$A_{t_1}=\left\{\widetilde{w}(s)\in\Gamma\right\}\bigcap
\left\{\eta_{t_1}\geq s\right\}.$ If $A_{t_1}=\emptyset$ then
$A_{t_1}\in{\Cal M}_{t_2}.$

Let $A_{t_1}\not=\emptyset.$ Then $A_{t_1}\bigcap\left\{\eta_{t_2}
\geq s\right\}\not=\emptyset,$ because $\{\eta_t\}$ is the
increasing process. Also $A_{t_1}=A_{t_2}\bigcap\{\eta_{t_1}\geq
s\}.$ Since $A_{t_2}\in\widetilde{{\Cal F}}_{t_2}$ and
$\{\eta_{t_1}\geq s\}\in{\Cal F}^w_{t_1} \subseteq{\Cal
F}^w_{t_2}$ then $A_{t_1}\in{\Cal M}_{t_2}.$

The lemma is proved.\enddemo

\proclaim{Lemma 2}The process $\{\zeta(t)\}$ is a square
integrable martingale w.r.t. $\{{\Cal M}_t\}$ and its
characteristic is equal to $\{\eta_t\}.$
\endproclaim
\demo{Proof} Firstly, we show that for all $t\geq 0$ the process
$\zeta(t)$ is ${\Cal M}_t$ measurable. We can approximate
$\zeta(t)$ by step functions in the following way
$$
\zeta(t)=\lim_{n\to\infty}\sum_{k=1}^{\infty}{\widetilde{w}\left({k-1\over
n} \right){\hbox{\rm 1\!I}}_{\left\{{k-1\over n}\leq\eta_t<{k\over
n}\right\}}}.
$$

We can write
$$\widetilde{w}\left({k-1\over n}\right){\hbox{\rm 1\!I}}_{\left\{{k-1\over
n}\leq\eta_t<{k\over
n}\right\}}=\left[\widetilde{w}\left({k-1\over n}\right){\hbox{\rm
1\!I}}_{\left\{{k-1\over
n}\leq\eta_t\right\}}\right]\left[{\hbox{\rm
1\!I}}_{\left\{\eta_t<{k\over
n}\right\}}\vphantom{\widetilde{w}\left({k-1\over n}\right)}
\right],$$ where the process in the first brackets is
$\widetilde{{\Cal F}}_t$ measurable and the process in the second
brackets is ${\Cal F}^w_t$ measurable. Then $\zeta(t)$ is ${\Cal
M}_t$ measurable as the limit of measurable functions.

The second moment of the process $\zeta(t)$ is finite for all
$t\geq 0$ because
$$
{\Bbb E}\widetilde{w}(\eta_t)^2={\Bbb E}{\Bbb
E}\left(\widetilde{w}(\eta_t)^2/{\Cal F}^w_\infty\right)={\Bbb
E}\left({\Bbb E}
\widetilde{w}(\widehat{t})^2\right)_{\widehat{t}=\eta_t}={\Bbb
E}\left(\widehat{t}\right)_{\widehat{t}= \eta_t}={\Bbb
E}\eta_t<\infty,
$$
in the second equality we used lemma 1, p.67, \cite{6}.

Let us prove that for arbitrary bounded ${\Cal M}_s$-measurable
random value $\xi$ the relation ${\Bbb
E}\left(\zeta(t)-\zeta(s)\right)\xi=0$ holds for all $t\geq s.$ It
is enough to check this relation for indicators of sets generating
$\sigma$-algebra ${\Cal M}_s.$ For all $k\geq 1,
u_i\in[0,+\infty),\Gamma_i\in{\Cal B}_S, 1\leq i\leq k,
\Gamma\in{\Cal B}_{C[0,+\infty)}$ (the Borel $\sigma$-algebra on
$C([0,+\infty))$)
$$
{\Bbb
E}\left[\widetilde{w}(\eta_t)-\widetilde{w}(\eta_s)\right]{\hbox{\rm
1\!I}}_{\left\{w(\cdot)|_0^s\in\Gamma
\right\}}\prod_{i=0}^{k}{\hbox{\rm
1\!I}}_{\left\{\widetilde{w}(u_i)\in\Gamma_i\right\}} {\hbox{\rm
1\!I}}_{\left\{u_i\leq\eta_s\right\}}=
$$
$$
={\Bbb E}\left[{\hbox{\rm
1\!I}}_{\left\{w(\cdot)|_0^s\in\Gamma\right\}}{\Bbb
E}\left(\left[\widetilde{w}
(\eta_t)-\widetilde{w}(\eta_s)\right]\prod_{i=0}^{k}{\hbox{\rm
1\!I}}_{\left\{\widetilde{w}(u_i)\in \Gamma_i\right\}}{\hbox{\rm
1\!I}}_{\left\{u_i\leq\eta_s\right\}}/{\Cal F}^w_{\infty}\right)
\right]=
$$
$$
={\Bbb E}\left.\left[{\hbox{\rm
1\!I}}_{\left\{w(\cdot)|_0^s\in\Gamma\right\}}{\Bbb
E}\left(\left[\widetilde{w}
(T_2)-\widetilde{w}(T_1)\right]\prod_{i=0}^{k}{\hbox{\rm
1\!I}}_{\left\{\widetilde{w}(u_i)\in\Gamma_i \right\}}{\hbox{\rm
1\!I}}_{\left\{u_i\leq
T_1\right\}}\right)\right]\right|_{T_2=\eta_t, T_1=\eta_s}=0,
$$
we denote by $w(\cdot)|_0^s$ the trajectory of the Wiener process
$\{w(t)\}$ on $[0,s].$ Here we use independence of the processes
$\{\widetilde{w}(t)\}, \{\eta_t\}$ and lemma 1, p.67, \cite{6} in
the second equality.

In the same way one can observe that the process
$\{\widetilde{w}(\eta_t)^2- \eta_t\}$ is martingale w.r.t.
$\{{\Cal M}_t\}.$ This means that the characteristics of the
martingale $\{\widetilde{w}(\eta_t)\}$ is equal to $\{\eta_t\}.$

Lemma is proved.\enddemo

\proclaim{Lemma 3}$\{w(t)\}$ is a Wiener process w.r.t. $\{{\Cal
M}_t\}.$\endproclaim

\remark{Remark 1} This result is not obvious because
$\sigma$-algebra $\{{\Cal M}_t\}$ is larger than $\{{\Cal
F}_t^w\}$ and $\sigma$-algebras $\{{\Cal F}_t^w\}$ and
$\{\widetilde{{\Cal F}}_t\}$ are not independent.
\endremark

\demo{Proof} It is enough to prove that the process $\{w(t)\}$ is
a martingale with characteristic $t$ w.r.t. $\{{\Cal M}_t\}.$
Firstly we show that for arbitrary bounded ${\Cal M}_s$ measurable
random value $\xi$ the relation ${\Bbb
E}\left(w(t)-w(s)\right)\xi=0$ holds for all $t\geq s.$ Again we
check this relation for indicators of sets generating
$\sigma$-algebra ${\Cal M}_s.$ For all $k\geq 1,
u_i\in[0,+\infty), \Gamma_i\in{\Cal B}_S, 1\leq i\leq k,
\Gamma\in{\Cal B}_{C[0,+\infty)}$
$$
{\Bbb E}\left[w(t)-w(s)\right]{\hbox{\rm
1\!I}}_{\left\{w(\cdot)|_0^s\in\Gamma\right\}}
\prod_{i=0}^{k}{\hbox{\rm
1\!I}}_{\left\{\widetilde{w}(u_i)\in\Gamma_i\right\}}{\hbox{\rm
1\!I}}_{\left\{u_i \leq\eta_s\right\}}=
$$
$$
={\Bbb E}\left[{\hbox{\rm
1\!I}}_{\left\{w(\cdot)|_0^s\in\Gamma\right\}}\prod_{i=0}^{k}{\hbox{\rm
1\!I}}_{ \left\{\widetilde{w}(u_i)\in\Gamma_i\right\}}{\hbox{\rm
1\!I}}_{\left\{u_i\leq\eta_s\right\}}{\Bbb E}
\left(w(t)-w(s)/{\Cal F}^w_s\vee{\Cal
F}^{\widetilde{w}}_{\infty}\right)\right]=0.
$$
In the same way one can show that the process $\{w(t)^2-t\}$ is
martingale w.r.t. $\{{\Cal M}_t\}.$  This means that the
characteristics of the martingale $\{w(t)\}$ is $t.$

Lemma is proved.\enddemo

For given $x_0\in\Re^d,$ measurable function $\alpha:S\to S$ and
operator $\beta:S\to{\Cal L}_+(S)$ (${\Cal L}_+(S)$ denotes the
space of all linear symmetric nonnegative operators on $S$) we
consider the following stochastic equation in $\Re^d$
$$
x(t)=x_0+\int_{0}^{t}{\left(q\nu+\alpha(x^S(\tau))\right)d\eta_\tau}+
\int_{0}^{t}{\widetilde{\beta}(x^S(\tau))d\widetilde{w}(\eta_\tau)}+w(t)\tag
1
$$
where $\widetilde{\beta}(\cdot)=\beta^{1/2}(\cdot).$ We call the
strong solution to the equation (1) the $\{{\Cal M}_t\}$-adapted
process $\{x(t)\}$ which satisfies the equality (1).

\proclaim{Theorem 1} Assume that there exists $K>0$ such that
\roster \item $ \sup_{x\in
S}\left(\left|\alpha(x)\vphantom{\widetilde{\beta}}\right|+\left\|\widetilde{\beta}(x)\right\|\right)\leq
K$, \item
$\left|\alpha(x)-\alpha(y)\vphantom{\widetilde{\beta}}\right|^2+\left\|\widetilde{\beta}(x)-\widetilde{\beta}(y)\right\|^2\leq
K\left|x-y\vphantom{\widetilde{\beta}}\right|^2,$ for all $x,y\in
S.$
\endroster
Then for all $T>0$ the solution to the equation (1) exists and is
unique when $t\in[0,T]$.
\endproclaim
It is enough to prove the existence and uniqueness for the process
$\{x^S(t)\}.$ This, in turn, is a consequence of the existence and
uniqueness theorems for the stochastic equations with arbitrary
martingales and stochastic measures (see \cite{7}, p. 278-296).

\remark{Remark 2} It was proved in \cite{8} that the solution of
the equation (1) has a measurable modification as the function of
the starting point.
\endremark

\head {3. The Markov property of the constructed process.}
\endhead

Let us denote by $\{x(t,x)\}$ the solution to the equation (1)
started from $x\in\Re^d.$

\proclaim{Theorem 2} $\{x(t,x)\}$ has Markov property.
\endproclaim

\demo{Proof} We prove the theorem if we show that the relation
$$
{\Bbb E}{\hbox{\rm 1\!I}}_{\{x(t,x)\in\Gamma\}}\xi={\Bbb
E}\Phi_{t-s}(x(s,x),\Gamma)\xi,\tag 2
$$
holds for an arbitrary bounded ${\Cal M}_s$-measurable random
value $\xi,$ where $\Phi_{t-s}(\cdot,\Gamma):\Re^d\to\Re$ is a
measurable function for all $0\leq s\leq t, \Gamma\in{\Cal
B}_{\Re^d}.$

At the beginning we note that the process $\{x(t,x)\}$ has the
property
$$
\theta_s x(t-s,z)|_{z=x(s,x)}=x(t,x),\quad 0\leq s\leq t,\tag 3
$$
here $\theta_s$ is the shift operator (see, \cite {10}, p.121). We
can deal with the process $\{x(t,x(s,x))\}$ because the process
$\{x(t,x)\}$ is measurable as the function of the initial point.
The equality (3) holds true because each of the processes at the
both sides of (3) satisfy the equation
$$
x(t,x)=x(s,x)+\int_s^t\left(q\nu+\alpha(x^S(u,x))\right)d\eta_u+
\int_s^t\widetilde{\beta}(x^S(u,x))
d\widetilde{w}(\eta_u)+w(t)-w(s),\tag 4
$$
that has the unique solution.

Let us denote
$(\gamma(\cdot)-\gamma(s))|_s^t=(\gamma(\cdot)-\gamma(s)){\hbox{\rm
1\!I}}_{[s,t]}(\cdot)$ for an arbitrary process $\{\gamma(t)\},
s\leq t.$ The equalities (3), (4) means that we can express the
process $\{x(t,x)\}$ in the form
$$
x(t,x)=F(x(s,x),(w(\cdot)-w(s))|_s^t,(\eta_\cdot-\eta_s)|_s^t,(\widetilde{w}(
\eta_\cdot)-\widetilde{w}(\eta_s))|_s^t),\quad 0\leq s\leq t,\tag
5
$$
where $F:\Re^d\times C[0,+\infty)\times C[0,+\infty)\to\Re^d$ is a
measurable functional. Therefore we prove (2) if we show that for
all $0\leq s\leq t,\Gamma\in{\Cal B}_{(C[0,+\infty))^3}$ the
relation
$$
{\Bbb E}{\hbox{\rm
1\!I}}_{\{\left((w(\cdot)-w(s))|_s^t,(\eta_\cdot-\eta_s)|_s^t,
(\widetilde{w}(\eta_\cdot)-\widetilde{w}(\eta_s))|_s^t\right)\in\Gamma\}}\xi=
{\Bbb E}\Phi^1_{t-s}(x(s,x),\Gamma)\xi, \tag 6
$$
holds true, where $\Phi^1_{t-s}(\cdot,\Gamma):\Re^d\to\Re$ is a
measurable function for all $t,s,\Gamma.$

We prove (6) in the next way. We construct
$\Phi^1_{t-s}(\cdot,\cdot)$ for the indicators of sets generating
$\sigma$-algebra ${\Cal M}_s.$ We show that for all  $k\geq 1,
u_i\in[0,+\infty),\widetilde{\Gamma}_i \in{\Cal B}_S,1\leq i\leq
k,\Gamma_0,\Gamma_1,\Gamma_2,\Gamma_3\in{\Cal B}_{C[0,+\infty)}$
the relation
$$
{\Bbb E}{\hbox{\rm
1\!I}}_{\{(w(\cdot)-w(s))|_s^t\in\Gamma_1\}}{\hbox{\rm
1\!I}}_{\{(\eta_\cdot-\eta_s)|_s^t \in\Gamma_2\}}{\hbox{\rm
1\!I}}_{\{(\widetilde{w}(\eta_\cdot)-\widetilde{w}(\eta_s))|_s^t\in\Gamma_3\}}
{\hbox{\rm
1\!I}}_{\{w(\cdot)|_0^s\in\Gamma_0\}}\prod_{i=1}^k\left({\hbox{\rm
1\!I}}_{\{\widetilde{w}(u_i)\in\widetilde
{\Gamma}_i\}}\right.\times
$$
$$
\times\left.{\hbox{\rm 1\!I}}_{\{u_i\leq\eta_s\}}\right)={\Bbb
E}\Phi^2_{t-s}(x(s,x),\Gamma_1, \Gamma_2,\Gamma_3){\hbox{\rm
1\!I}}_{\{w(\cdot)|_0^s\in\Gamma_0\}}\prod_{i=1}^k{\hbox{\rm
1\!I}}_{\{\widetilde{w}(u_i) \in\widetilde{\Gamma}_i\}}{\hbox{\rm
1\!I}}_{\{u_i\leq \eta_s\}}, \tag 7
$$
holds, where
$\Phi^2_{t-s}(\cdot,\Gamma_1,\Gamma_2,\Gamma_3):\Re^d\to\Re$ is a
measurable function for all $t,s,\Gamma_1,\Gamma_2,\Gamma_3.$

Let us take the conditional expectation w.r.t. ${\Cal F}_s^w$ on
the left hand side of the relation (7). Then we obtain
$$
{\Bbb E}{\hbox{\rm
1\!I}}_{\{(w(\cdot)-w(s))|_s^t\in\Gamma_1\}}{\hbox{\rm
1\!I}}_{\{(\eta_\cdot-\eta_s)|_s^t \in\Gamma_2\}}{\hbox{\rm
1\!I}}_{\{(\widetilde{w}(\eta_\cdot)-\widetilde{w}(\eta_s))|_s^t\in\Gamma_3\}}
{\hbox{\rm 1\!I}}_{\{w(\cdot)|_0^s\in\Gamma_0\}}\times
$$
$$
\times\prod_{i=1}^k{\hbox{\rm
1\!I}}_{\{\widetilde{w}(u_i)\in\widetilde{\Gamma}_i\}}{\hbox{\rm
1\!I}}_{\{u_i\leq\eta_s \}}={\Bbb E}\left({\hbox{\rm
1\!I}}_{\{w(\cdot)|_0^s\in\Gamma_0\}}{\Bbb E}\left[{\hbox{\rm
1\!I}}_{\{(w(\cdot)-z) |_s^t\in\Gamma_1\}}{\hbox{\rm
1\!I}}_{\{(\eta_\cdot-r)|_s^t\in\Gamma_2\}}\right.\right. \times
$$
$$
\times\left.\left.\left.{\hbox{\rm
1\!I}}_{\{(\widetilde{w}(\eta_\cdot)-\widetilde{w}(r))|_s^t\in
\Gamma_3\}}\prod_{i=1}^k{\hbox{\rm
1\!I}}_{\{\widetilde{w}(u_i)\in\widetilde{\Gamma}_i\}}{\hbox{\rm
1\!I}}_{\{u_i\leq r\}}/ {\Cal F}_s^w
\right]\right|_{z=w(s),r=\eta_s}\right)\tag 8
$$

Let us show that the sets
$\{(w(\cdot)-z)|_s^t\in\Gamma_1\}$$\bigcap\{(
\eta_\cdot-r)|_s^t\in\Gamma_2\}$$\bigcap\{(\widetilde{w}(\eta_\cdot)-\widetilde{w}(r))
|_s^t\in\Gamma_3\}$ and
$\bigcap_{i=1}^k{\{\widetilde{w}(u_i)\in\widetilde{\Gamma}_i\}}
\bigcap{\{u_i\leq r\}}$ are independent when the $\sigma$-algebra
${\Cal F}_s^w$ is fixed. Really, for arbitrary bounded ${\Cal
F}_s^w$-measurable random value $\zeta$ we have
$$
{\Bbb E}\zeta{\hbox{\rm
1\!I}}_{\{(w(\cdot)-z)|_s^t\in\Gamma_1\}}{\hbox{\rm
1\!I}}_{\{(\eta_\cdot-r)|_s^t\in\Gamma_2 \}}{\hbox{\rm
1\!I}}_{\{(\widetilde{w}(\eta_\cdot)-\widetilde{w}(r))|_s^t\in\Gamma_3\}}\prod_{i=1}^k
{\hbox{\rm
1\!I}}_{\{\widetilde{w}(u_i)\in\widetilde{\Gamma}_i\}}{\hbox{\rm
1\!I}}_{\{u_i\leq r\}}=
$$
$$
={\Bbb E}\left(\zeta{\hbox{\rm
1\!I}}_{\{(w(\cdot)-z)|_s^t\in\Gamma_1\}}{\hbox{\rm
1\!I}}_{\{(\eta_\cdot-r)|_s^t \in\Gamma_2\}}{\Bbb
E}\left[{\hbox{\rm
1\!I}}_{\{(\widetilde{w}(\eta_\cdot)-\widetilde{w}(r))|_s^t\in\Gamma_3\}}
\right.\right.\times
$$
$$
\times\left.\left.\prod_{i=1}^k{\hbox{\rm
1\!I}}_{\{\widetilde{w}(u_i)\in\widetilde{\Gamma}_i\}}{\hbox{\rm
1\!I}}_{\{ u_i\leq r\}}/{\Cal F}_\infty^w\right]\right)={\Bbb
E}\left(\zeta{\hbox{\rm 1\!I}}_{\{(w(\cdot)-z)
|_s^t\in\Gamma_1\}}{\hbox{\rm
1\!I}}_{\{(\eta_\cdot-r)|_s^t\in\Gamma_2\}}\right.\times
$$
$$
\times\left.\left.{\Bbb E}\left[{\hbox{\rm
1\!I}}_{\{(\widetilde{w}(r(\cdot))-\widetilde{w}(r))|_s^t\in
\Gamma_3\}}\prod_{i=1}^k{\hbox{\rm
1\!I}}_{\{\widetilde{w}(u_i)\in\widetilde{\Gamma}_i\}}{\hbox{\rm
1\!I}}_{\{u_i\leq r\}}
\right]\right|_{r(\cdot)=\eta_\cdot}\right)={\Bbb
E}\left(\zeta{\hbox{\rm 1\!I}}_{\{(w(\cdot)
-z)|_s^t\in\Gamma_1\}}\right.\times
$$
$$
\times\left.\left.{\hbox{\rm
1\!I}}_{\{(\eta_\cdot-r)|_s^t\in\Gamma_2\}}{\Bbb
E}\left[{\hbox{\rm 1\!I}}_{\{(
\widetilde{w}(r(\cdot))-\widetilde{w}(r))|_s^t\in\Gamma_3\}}\right|_{r(\cdot)=\eta_\cdot}
\right]{\Bbb E}\left[\prod_{i=1}^k{\hbox{\rm
1\!I}}_{\{\widetilde{w}(u_i)\in\widetilde{\Gamma}_i\}}{\hbox{\rm
1\!I}}_{\{u_i \leq r\}}\right]\right)=
$$
$$
={\Bbb E}\left(\zeta{\hbox{\rm
1\!I}}_{\{(w(\cdot)-z)|_s^t\in\Gamma_1\}}{\hbox{\rm
1\!I}}_{\{(\eta_\cdot-r)|_s^t \in\Gamma_2\}}{\Bbb
E}\left[{\hbox{\rm
1\!I}}_{\{(\widetilde{w}(\eta_\cdot)-\widetilde{w}(r))|_s^t\in\Gamma_3\}}/
{\Cal F}_\infty^w\right]\right)\times
$$
$$
\times{\Bbb E}\left(\prod_{i=1}^k{\hbox{\rm
1\!I}}_{\{\widetilde{w}(u_i)\in\widetilde{\Gamma}_i\}}{\hbox{\rm
1\!I}}_{\{u_i \leq r\}}\right)={\Bbb E}\left(\zeta{\hbox{\rm
1\!I}}_{\{(w(\cdot)-z)|_s^t\in\Gamma_1\}}{\hbox{\rm 1\!I}}_{\{(
\eta_\cdot-r)|_s^t\in\Gamma_2\}}\right.\times
$$
$$
\times\left.{\hbox{\rm
1\!I}}_{\{(\widetilde{w}(\eta_\cdot)-\widetilde{w}(r))|_s^t\in\Gamma_3\}}\right)
{\Bbb E}\left(\prod_{i=1}^k{\hbox{\rm
1\!I}}_{\{\widetilde{w}(u_i)\in\widetilde{\Gamma}_i\}}{\hbox{\rm
1\!I}}_{\{u_i\leq r\} }\right).
$$
In the third equality we take into account independence of the
sets
$\{(\widetilde{w}(r(\cdot))-\widetilde{w}(r))|_s^t\in\Gamma_3\}$
and
$\bigcap_{i=1}^k{\{\widetilde{w}(u_i)\in\widetilde{\Gamma}_i\}}\bigcap{\{u_i\leq
r\}}.$

Therefore we obtain that the expression on the right-hand side of
(8) equals
$$
{\Bbb E}\left({\hbox{\rm
1\!I}}_{\{w(\cdot)|_0^s\in\Gamma_0\}}{\Bbb E}\left[{\hbox{\rm
1\!I}}_{\{(w(\cdot)-w(s)) |_s^t\in\Gamma_1\}}{\hbox{\rm
1\!I}}_{\{(\eta_\cdot-\eta_s)|_s^t\in\Gamma_2\}}{\hbox{\rm
1\!I}}_{\{(\widetilde{w}
(\eta_\cdot)-\widetilde{w}(\eta_s))|_s^t\in\Gamma_3\}}/{\Cal
F}_s^w\right]\right.\times
$$
$$
\times\left.{\Bbb E}\left[\prod_{i=1}^k{\hbox{\rm
1\!I}}_{\{\widetilde{w}(u_i)\in\widetilde{\Gamma}_i\}}{\hbox{\rm
1\!I}}_{\{ u_i\leq\eta_s\}}/{\Cal F}_s^w \right]\right).\tag 9
$$

Let us denote the conditional joint distribution of the processes
$\{(w(\cdot)-w(s))|_s^t\}$ and $\{(\eta_\cdot-\eta_s)|_s^t\}$ in
the following way
$$
{\Bbb
P}\left\{(w(\cdot)-w(s))|_s^t\in\Gamma_1,(\eta_\cdot-\eta_s)|_s^t\in\Gamma_2/
{\Cal F}_s^w\right\}=\mu_{t-s}(x^\nu(s),\Gamma_1,\Gamma_2).
$$
We remind that the process $\{\eta_t\}$ is an additive functional
of the $\{x^\nu (t)\}$ and thus the distribution of
$\{(\eta_\cdot-\eta_s)|_s^t\}$ depends only on distribution of
$x^\nu(s).$ Then we have
$$
{\Bbb E}\left({\hbox{\rm
1\!I}}_{\{(w(\cdot)-w(s))|_s^t\in\Gamma_1\}}{\hbox{\rm
1\!I}}_{\{(\eta_\cdot-\eta_s) |_s^t\in\Gamma_2\}}{\hbox{\rm
1\!I}}_{\{(\widetilde{w}(\eta_\cdot)-\widetilde{w}(\eta_s))|_s^t\in\Gamma_3\}}
/{\Cal F}_s^w\right)=
$$
$$
=\int{\Bbb
P}\left\{y(\cdot)\in\Gamma_1,\theta(\cdot)\in\Gamma_2,\widehat{w}(\theta(\cdot))\in
\Gamma_3\right\}\mu_{t-s}(x^\nu(s),dy(\cdot),d\theta(\cdot))=
$$
$$
=\Phi^2_{t-s}(x^\nu(s),\Gamma_1,\Gamma_2,\Gamma_3),
$$
where $\widehat{w}(\cdot)$ is a Wiener process in $S$.

Therefore (9) equals
$$
{\Bbb E}\left({\hbox{\rm
1\!I}}_{\{w(\cdot)|_0^s\in\Gamma_0\}}{\Bbb E}\left[{\hbox{\rm
1\!I}}_{\{(w(\cdot)-w(s)) |_s^t\in\Gamma_1\}}{\hbox{\rm
1\!I}}_{\{(\eta_\cdot-\eta_s)|_s^t\in\Gamma_2\}}{\hbox{\rm
1\!I}}_{\{(\widetilde{w}
(\eta_\cdot)-\widetilde{w}(\eta_s))|_s^t\in\Gamma_3\}}/{\Cal
F}_s^w\right]\right.\times
$$
$$
\times\left.{\Bbb E}\left[\prod_{i=1}^k{\hbox{\rm
1\!I}}_{\{\widetilde{w}(u_i)\in\widetilde{\Gamma}_i\}}{\hbox{\rm
1\!I}}_{\{ u_i\leq\eta_s\}}/{\Cal F}_s^w\right]\right)={\Bbb
E}\left({\hbox{\rm 1\!I}}_{\{w(\cdot)|_0^s\in
\Gamma_0\}}\Phi^2_{t-s}(x^\nu(s),\Gamma_1,\Gamma_2,\Gamma_3)\right.\times
$$
$$
\times\left.{\Bbb E}\left[\prod_{i=1}^k{\hbox{\rm
1\!I}}_{\{\widetilde{w}(u_i)\in\widetilde{\Gamma}_i\}}{\hbox{\rm
1\!I}}_{\{ u_i\leq\eta_s\}}/{\Cal F}_s^w\right]\right)={\Bbb
E}\left({\hbox{\rm 1\!I}}_{\{w(\cdot)|_0^s\in
\Gamma_0\}}\Phi^2_{t-s}(x^\nu(s),\Gamma_1,\Gamma_2,\Gamma_3)\right.\times
$$
$$
\times\left.\prod_{i=1}^k{\hbox{\rm
1\!I}}_{\{\widetilde{w}(u_i)\in\widetilde{\Gamma}_i\}}{\hbox{\rm
1\!I}}_{\{u_i\leq \eta_s\}}\right).
$$
This proves the equality (7) and then, after standard limiting
procedure, the equality (6). The theorem is proved.
\enddemo

In the next theorem we characterize the process constructed in the
the Theorem 1 as a generalized diffusion process in the sense of
Portenko (see \cite{1}).

\proclaim{Theorem 3} For all $\theta\in\Re^d, \varphi\in{\Cal
C}_0(\Re^d)$ (the space of all real-valued continuous functions in
$\Re^d$ with compact support) the relations
$$
\lim_{t\downarrow 0}{1\over t}\int_{\Re^d}\varphi(z) {\Bbb
E}|x(t,z)-x(0,z)|^4dz=0,
$$
$$
\lim_{t\downarrow 0}{1\over t}\int_{\Re^d}{\varphi(z){\Bbb
E}(x(t,z)- x(0,z),\theta)dz}=\int_{S}{\varphi(z)(q\nu
+\alpha(z),\theta)d\sigma_z}
$$
$$
\lim_{t\downarrow 0}{1\over t}\int_{\Re^d}{\varphi(z) {\Bbb
E}(x(t,z)-x(0,z),\theta)^2dz}=(\theta,\theta)\int_{\Re^d}{\varphi(z)dz}+\int_{S}{(\beta(z)\theta,\theta)
\varphi(z)d\sigma_z}.
$$
hold, where $d\sigma$ is the Lebesgue measure on $S$.\endproclaim

This result means that the process $\{x(t,x)\}$ is a generalized
diffusion process with the generalized drift vector being equal to
$(q\nu+\alpha(x^S))\delta_S(x)$ and the generalized diffusion
matrix being equal to $I+\beta(x^S)\delta_S(x),$ here $I$ is the
identity matrix in $\Re^d,$ $\delta_S(\cdot)$ is the
delta-function concentrated on $S$.

The proof of this theorem is similar to the corresponding theorem
in \cite{11}.

\end